\documentclass[11pt]{article}

\usepackage[cp850]{inputenc}

\usepackage{amscd}      
\usepackage{amssymb}
\usepackage[intlimits]{amsmath}

\usepackage{theorem}

\newtheorem{prop}{Proposition}[section]
\newtheorem{lem}[prop]{Lemma}

\newtheorem{them}[prop]{Theorem}

\theorembodyfont{\upshape}

\newtheorem{defn}[prop]{Definition}

\newtheorem{numrmk}[prop]{Remark}

\newtheorem{numex}[prop]{Example}

\newenvironment{pf}{\begin{trivlist}\item[]{\sc Proof.}}%
           {\nolinebreak $\Box$ \end{trivlist}}

\newcommand{\noprint}[1]{}

\newcommand{\frakd}{{\mathfrak d}}
\newcommand{\frakg}{{\mathfrak g}}

\newcommand{\G}{{\cal G}}

\newcommand{\hH}{{\cal H}}

\newcommand{\rR}{{\cal R}}
\newcommand{\TX}{{\mathbb{T}X}}
\newcommand{\TS}{{\mathbb{T}S}}
\newcommand{\TM}{{\mathbb{T}M}}
\newcommand{\TG}{{\mathbb{T}\mathcal{G}}}

\newcommand{\M}{{\overline{\mathcal{M}}}}

\newcommand{\T}{{\mathbb{T}}}
\newcommand{\pr}{{\mathrm{pr}}}

\newcommand{\Cour}[1]      {[\![#1]\!]}

\def\gpd{\rightrightarrows}

\newcommand{\s}{{\bf s}}             
\renewcommand{\t}{{\bf t}}           

\newcommand{\Lie}{\mathcal{L}}

\newcommand{\gm}{\Gamma }

\newcommand{\be}{\begin{eqnarray*}}
\newcommand{\ee}{\end{eqnarray*}}

\newcommand{\D}{K}  

\newcommand{\SP} [1]     {{\left\langle {{#1}} \right\rangle}}

\begin{document}
\title{{\bf Courant morphisms and moment maps}}
\author{Henrique Bursztyn$^1$, David Iglesias Ponte$^2$, Pavol
\v{S}evera$^3$
\\{\small\it $^1$Instituto Nacional de Matem\'atica Pura e Aplicada (IMPA),
Brazil}
\\{\small\it $^2$Instituto de Matem\'aticas y F\'{\i}sica Fundamental (CSIC),
Spain}
\\{\small\it $^3$Section de Math\'ematiques, Universit\'e de Gen\`eve, Switzerland,}
\\{\small\it on leave from FMFI UK, Bratislava, Slovakia}
\\[5pt]
{\small\it e-mail: henrique@impa.br, iglesias@imaff.cfmac.csic.es,
pavol.severa@gmail.com} }
\date{}

\sloppy
\maketitle

\begin{abstract}
We study Hamiltonian spaces associated with pairs $(E,A)$, where
$E$ is a Courant algebroid and $A\subset E$ is a Dirac structure.
These spaces are defined in terms of morphisms of Courant
algebroids with suitable compatibility conditions. Several of
their properties are discussed, including a reduction procedure.
This set-up encompasses familiar moment map theories, such as
group-valued moment maps, and it provides an intrinsic approach
from which different geometrical descriptions of moment maps can
be naturally derived. As an application, we discuss the
relationship between quasi-Poisson and presymplectic groupoids.
\end{abstract}

\tableofcontents

\section{Introduction}
In this paper, we study Hamiltonian spaces naturally associated
with \emph{Manin pairs}, i.e.\ pairs $(E,A)$, where $E$ is a Courant algebroid, and
$A\subset E$ is a Dirac structure. As we will see, these objects
retain many of the features of ordinary Hamiltonian spaces,
including a reduction procedure, and several moment map theories
can be expressed in these terms.

Our main motivation to consider this general set-up comes from the
theory of Hamiltonian spaces with group-valued moment maps (as in
\cite{Lu} and \cite{AKM,AMM}) or, more generally, of moment maps
with values in homogeneous spaces $D/G$. There are two geometrical
formulations of $D/G$-valued moment maps: The original one of
Alekseev and Kosmann-Schwarzbach \cite{AK} is given in terms of
\textit{quasi-Poisson actions}, whereas \cite{BC2} presents an
alternative approach based on morphisms of \textit{Dirac
manifolds}. (In the special case of $G$-valued moment maps, these
two distinct viewpoints can be found in \cite{AKM} and \cite{AMM},
respectively.) Both approaches originate from a common starting
point, namely a \textit{Manin pair}, but each one requires an
additional, noncanonical piece of information; depending on this
extra choice, Hamiltonian spaces are either described in terms of
quasi-Poisson or twisted Dirac structures. As proven in
\cite[Sec.~6]{BC2}, these two viewpoints, despite resorting to
different geometrical structures, produce isomorphic categories of
Hamiltonian spaces, regardless of any of the extra data used. This
raises the question of whether there is an \textit{intrinsic}
notion of Hamiltonian space associated with a Manin pair,
requiring no extra data at all, that would naturally recover the
formulations in \cite{AK} and \cite{BC2} once suitable additional
choices are made. The main goal of this paper is to present such
intrinsic notion and to study its properties, showing that it
offers a clear conceptual explanation for the equivalence between
the quasi-Poisson and Dirac geometric approaches to moment maps.

We organize the paper as follows. We review the basics of Courant
algebroids and Dirac structures \cite{LWX} in Section
\ref{sec:prelim}, including Dirac structures supported on a
submanifold and Courant algebroid morphisms \cite{AX}, and we
introduce the notion of \textit{morphism of Manin pairs}, which
plays a central role in this paper. In Section \ref{sec:ham}, we
define \textit{Hamiltonian spaces} associated with Manin pairs
$(E,A)$ over a manifold $S$; these objects are special examples of
morphisms of Manin pairs. More explicitly, Hamiltonian spaces are
triples $(X,J,\D)$, where $X$ is a manifold, $J:X\to S$ is a
smooth map (the \textit{moment map}), and $\D$ is a Dirac
structure on the product Courant algebroid $(TX\oplus T^*X) \times
E$ with support on $\mathrm{graph}(J)\subset X\times S$,
satisfying suitable compatibility conditions. We focus on two
possible scenarios: when $E=A\oplus A^*$ is the double of a Lie
quasi-bialgebroid \cite{Roy}, then the Hamiltonian spaces for
$(E,A)$ can be naturally identified with Hamiltonian quasi-Poisson
spaces (in the sense of \cite{ILX}); on the other hand, when
$E=TS\oplus T^*S$ is a Courant algebroid defined by a closed
3-form on $S$ \cite{SeWe01}, then $A\subset TS\oplus T^*S$ is a
Dirac structure on $S$, and the Hamiltonian spaces for $(E,A)$ are
identified with Hamiltonian spaces defined by strong Dirac maps
into $S$ (as considered in \cite{ABM,BC,BC2}). Combining these two
possible ``incarnations'' of the intrinsic Hamiltonian spaces for
$(E,A)$, we arrive at a functorial correspondence between moment
maps in quasi-Poisson and Dirac geometries, which gives a new,
more general viewpoint to the results in \cite[Sec.~6.2]{BC2} (as
well as \cite[Sec.~5.4]{ABM} and \cite[Sec.~3.5]{BC}). We also
discuss Poisson spaces obtained via reduction of Hamiltonian
spaces for $(E,A)$, showing that they agree with quasi-Poisson and
Dirac reductions in specific situations. Finally, in Section
\ref{sec:integ}, we apply the correspondence between Hamiltonian
spaces in quasi-Poisson and Dirac geometries to find an explicit
construction relating presymplectic groupoids \cite{BCWZ} and
quasi-Poisson groupoids \cite{ILX}.

\medskip

\noindent{\bf Acknowledgments:} We thank P. Xu for stimulating
discussions and Zhuo Chen and the referee for their comments. We acknowledge
the financial support of CNPq (H.B.), MEC Research Contract ``Juan
de la Cierva", grant MTM2007-62478 (D.I.P.) and the Swiss National
Science Fundation (P.\v S.). We also thank several institutions
for hosting us while this work was being done, including IMPA
(D.I.P.), University of La Laguna (H.B. and D.I.P.) and the Erwin
Schr\"odinger Institute.


\section{Preliminaries}\label{sec:prelim}

\subsection{Courant algebroids and Manin pairs}

A \textbf{Courant algebroid} \cite{LWX} over a manifold $S$ is a
vector bundle $E\to S$ equipped with a nondegenerate symmetric
bilinear form $\SP{ \cdot , \cdot}$ on the bundle, a bundle map
$\rho :E\to TS$ and a bilinear bracket $\Cour{ \cdot, \cdot }$ on
$\Gamma (E)$ such that $ \forall \,e, e_1,e_2,e_3\in \Gamma (E)$,
$f\in C^\infty (S)$, the following holds:
\begin{itemize}
\item[$c1)$] $\Cour{ e_1,\Cour{ e_2,e_3} } =\Cour{ \Cour{ e_1,
e_2},e_3 } +\Cour{ e_2,\Cour{ e_1,e_3} }$;

\item[$c2)$] $\Cour{ e,e} = \rho^*d\SP{ e,e}$, where we use
$\SP{\cdot,\cdot}$ to identify $E\cong E^*$;

\item[$c3)$] ${\cal L}_{\rho (e)}\SP{ e_1,e_2} =\SP{ \Cour{
e,e_1},e_2}+\SP{ e_1,\Cour{ e,e_2} }$;

\item[$c4)$] $\rho (\Cour{ e_1,e_2} )=[ \rho (e_1), \rho (e_2)]$;

\item[$c5)$] $\Cour{ e_1,f\, e_2} =f\Cour{ e_1, e_2} + ({\cal
L}_{\rho (e_1)}f)e_2$.
\end{itemize}
A Courant algebroid is denoted by the quadruple
$(E,\SP{\cdot,\cdot},\Cour{\cdot,\cdot},\rho)$, or simply by $E$
if there is no risk of confusion. We use the notation
$\overline{E}$ for the Courant algebroid
$(E,-\SP{\cdot,\cdot},\Cour{\cdot,\cdot},\rho)$.

We recall some properties of Courant algebroids for later use.
First, $c2)$ and $c4)$ imply that
\begin{equation}\label{eq:rhorho*}
\rho\circ \rho^*=0,
\end{equation}
where we identify $E\cong E^*$ via $\SP{\cdot,\cdot}$. On the
other hand, given 1-forms $\beta, \beta' \in \Omega^1(S)$, we have
\begin{equation}\label{eq:prop2}
\Cour{\rho^*\beta,\rho^*\beta'}=0.
\end{equation}
To verify \eqref{eq:prop2}, it suffices to show that
$\SP{\Cour{\rho^*\beta,\rho^*\beta'},e}=0, \, \forall \, e\in
\Gamma(E)$, and this follows from $c3)$ combined with $c4)$ and
\eqref{eq:rhorho*}.

A subbundle $A\subset E$ is called an \textbf{almost Dirac
structure} if it is Lagrangian with respect to $\langle
\cdot ,\cdot \rangle$ (i.e.\ both isotropic and coisotropic),
and it is a \textbf{Dirac structure} if, in
addition, it is \textit{integrable}, that is, $\Gamma (A)$ is
closed under $\Cour{ \cdot ,\cdot}$. The restrictions of the
anchor and Courant bracket to any Dirac structure $A$ make it into
a \textit{Lie algebroid}. We denote the restricted bracket by
$[\cdot,\cdot]_A:=\Cour{\cdot,\cdot}|_{\Gamma(A)}$.
 Pairs $(E,A)$, where $E$ is a Courant
algebroid over $S$ and $A\subset E$ is a Dirac structure, are
central objects in this paper. We refer to them as \textbf{Manin
pairs over $S$}.\footnote{Lagrangian subbundles $A\subset E$ exist
if and only if the pairing $\SP{\cdot,\cdot}$ has split signature
$(n,n)$. The results in this paper remain valid when $A\subset E$
is only required to be maximal isotropic, in which case there is
no signature requirement on $\SP{\cdot,\cdot}$.}

\begin{numex}\label{ex:point}
A Manin pair over a point is a pair $(\frakd,\frakg)$ \cite{AK},
where $\frakd$ is a $2n$-dimensional Lie algebra equipped with an
Ad-invariant inner product, and
$\frakg\subset \frakd$ is a Lagrangian Lie subalgebra.
\end{numex}

\begin{numex}\label{ex:dressing}
Let $(\frakd,\frakg)$ be a Manin pair (over a point) as in
Example~\ref{ex:point}. Following \cite{AK}, let $D$ and $G$ be
connected Lie groups integrating $\frakd$ and $\frakg$, and assume
that $G$ is a closed subgroup of $D$. Let $S:=D/G$ be the quotient
with respect to right multiplication. Then the action of $D$ on
itself by left multiplication induces a $D$-action on $S$, called
the \textbf{dressing action}. The pairing and bracket on $\frakd$
give rise to a natural Courant algebroid structure on the trivial
bundle $\frakd_S:=\frakd\times S$ over $S$ \cite{AX,Se} (c.f.
\cite[Sec.~3]{BC2}) for which the anchor is the infinitesimal
dressing action $\frakd_S\to TS$. The subbundle
$\frakg_S=\frakg\times S$ defines a Dirac structure, i.e.\
$(\frakd_S,\frakg_S)$ is a Manin pair over $S$.
\end{numex}

\begin{numex}\label{ex:standard}
Let $\TS:=TS\oplus T^*S$, equipped with pairing
$\SP{(v,\alpha),(v',\alpha')}=\alpha'(v)+\alpha(v')$. A closed
3-form $\phi_S\in \Omega^3(S)$ defines the Courant bracket
\cite{SeWe01}
$$
\Cour{(v,\alpha),(v',\alpha')}:=([v,v'],\Lie_v\alpha'-i_{v'}{d\alpha}+i_{v'}i_v\phi_S),
$$
making $\TS$ a Courant algebroid with anchor given by the natural
projection $\TS\to TS$. When $\phi_S=0$, we refer to this Courant
structure on $\TS$ as \textbf{standard}.
\end{numex}

Let $(E,A)$ be a Manin pair over $S$. Consider the exact sequence
$0\to A \to E\to E/A\cong A^*\to 0$. We can always fix an
\textit{isotropic splitting} of this sequence, i.e., a splitting
\begin{equation}\label{eq:j}
j:A^* \to E
\end{equation}
whose image is isotropic (see e.g. \cite[App.~2]{BC2}). This
defines a cobracket $F_j:\Gamma(A)\to \wedge^2\Gamma(A)$, a bundle
map $\rho_{A^*}^j:A^*\to TM$ and a 3-tensor $\chi_j\in
\Gamma(\wedge^3A)$ by
\begin{equation}\label{eq:liequasi}
F_j(a)(\xi_1,\xi_2)=\SP{\Cour{j(\xi_1),j(\xi_2)},a},\;\;
\rho_{A^*}^j=\rho\circ j,\;\;
\chi_j(\xi_1,\xi_2,\xi_3)=\SP{\Cour{j(\xi_1),j(\xi_2)},j(\xi_3)},
\end{equation}
where $a\in \Gamma(A),\, \xi_1, \xi_2, \xi_3\in \Gamma(A^*)$,
making $A$ into a \textbf{Lie quasi-bialgebroid} \cite{ILX,Roy} in
such a way that $E$ is naturally identified with the ``double''
Courant algebroid $A\oplus A^*$ \cite{Roy}. Equivalently, $F_j$
and $\rho_{A^*}^j$ can be combined into a degree 1 derivation
$d_{A^*}:\Gamma(\wedge^\bullet A)\to \Gamma(\wedge^{\bullet +1}
A)$ such that $d_{A^*}^2=[\chi_j,\cdot]_A$ and $d_{A^*}\chi_j=0$.

\subsection{Dirac structures with support}\label{subsec:diracsupp}

\begin{defn}\label{def:dirac}Given a Courant algebroid
$(E,\SP{\cdot,\cdot},\Cour{\cdot,\cdot},\rho)$ over a manifold $M$
and a submanifold $\iota: Q \hookrightarrow M$, a \textbf{Dirac
structure supported on $Q$} (see \cite{AX,Some}) is a subbundle $\D
\subset \iota^*E=E|_Q$ such that $\D_x\subset E_x$ is maximal
isotropic for all $x\in Q$ and the following two conditions hold:
\begin{itemize}
\item[$d1)$] $\D$ is compatible with the anchor, that is,
$\rho(\D)\subset TQ$;

\item[$ d2)$] for any sections $e_1,e_2$ of $E$ such that
$e_1|_Q$, $e_2|_Q\in \gm (\D)$, $\Cour{ e_1,e_2} |_Q\in \gm (\D)$.
\end{itemize}
If only $d1)$ is satisfied, we refer to an \textbf{almost Dirac
structure supported on $Q$}.
\end{defn}


\begin{numrmk}\label{rmk:bases}
The Leibniz rule for Courant algebroids (condition $c5)$) shows
that if $e_1$ and $e_1'$ satisfy $e_1|_Q=e_1'|_Q$, then $\Cour{e,e_1} |_Q= \Cour{e, e_1'} |_Q$ for all $e\in \Gamma(E)$ with $e|_Q$ in $\D$.
Hence it suffices to check $d2)$ for a set of sections of $E$
whose restrictions to $Q$ locally generate $K$.
\end{numrmk}

\begin{numex}
Let $Q$ be a submanifold of $M$, and denote by $NQ$ its normal
bundle. We consider $\TM$ equipped with its standard Courant
algebroid structure. Then the subbundle $TQ\oplus NQ \subset
\TM|_Q$ is a Dirac structure supported on $Q$.
\end{numex}

\begin{numex}\label{ex:gammaf}
Let $X$ and $X'$ be smooth manifolds, and let $f:X\to X'$ be a
smooth map. We consider the standard Courant algebroids $\TX'$ and
$\overline{\TX}$, and the product Courant algebroid $E=\TX' \times
\overline{\TX}$ over $M=X'\times X$. Then
\begin{equation}\label{gds-basic}
\Gamma_f:=\{ ((df(v),\alpha ),(v,f ^*\alpha)) \,|\,  v\in T_xX
\mbox{ and }\alpha \in T_{f(x)}^*X',\, x\in X \}
\end{equation}
is a Dirac structure in $E$ supported on
$Q=\mathrm{graph}(f)\subset M$.
\end{numex}

If $E$, $E'$ are Courant algebroids over $X$, $X'$, then a
\textbf{Courant algebroid morphism} between $E$ and $E'$ is a
Dirac structure in $E\times \overline{E}'$ supported on
$\mathrm{graph}(f)$, where $f:X\to X'$ is a smooth map (see e.g.
\cite{AX} and Remark \ref{rem:super}).

We have the following natural correspondence.

\begin{prop}\label{prop:correspondence}
Let $f:Q \hookrightarrow M$ be an embedding and $\phi\in
\Omega^3(M)$ be a closed 3-form. Let $\mathrm{Dir}(Q)$ denote the
set of Dirac structures in $\mathbb{T}Q$, integrable relative to
$f^*\phi$, and $\mathrm{Dir}(M)_{f(Q)}$ denote the set of Dirac
structures in $\TM$ supported on $f(Q)$, integrable relative to
$\phi$. Then these sets are in bijection via the maps
\begin{align*}
&\mathfrak{F}_f:\mathrm{Dir}(Q)\to \mathrm{Dir}(M)_{f(Q)},\;\;
\mathfrak{F}_f(L)_{f(x)}:=\{ (df(u),\beta)\;|\; (u,f^*\beta)\in L_x \},\\
&\mathfrak{B}_f:\mathrm{Dir}(M)_{f(Q)}\to \mathrm{Dir}(Q),\;\;
\mathfrak{B}_f(L)_{x}:=\{ (u,f^*\beta)\;|\; (df(u),\beta)\in
L_{f(x)}\},
\end{align*}
which are inverses of one another.
\end{prop}

\begin{pf}
Consider the injective bundle map $\psi:\mathbb{T}Q\to f^*TM\oplus
T^*Q$, $\psi(u,\alpha)=(df(u),\alpha)$, and the surjective bundle
map $p: f^*\mathbb{T}M\to f^*TM\oplus T^*Q$,
$p(v,\beta)=(v,f^*\beta)$. Then
$\mathfrak{F}_f(L)=p^{-1}(\psi(L))$ is a smooth vector bundle over
$f(Q)$, which one can directly check to be maximal isotropic. To
verify integrability, consider sections $(v,\beta), (v',\beta')$
of $\mathbb{T}M$ whose restrictions to $f(Q)$ lie in
$\mathfrak{F}_f(L)$. Let $u, u'\in \Gamma(TQ)$ be such that
$v|_{f(Q)}=df(u), v'|_{f(Q)}=df(u')$, and consider
$\Cour{(v,\beta),(v',\beta')}=([v,v'],\mathcal{L}_v\beta'-i_{v'}d\beta
+i_{v\wedge v'}\phi)$. Then $[v,v']|_{f(Q)}=df([u,u'])\in Tf(Q)$
and, similarly, $f^*\mathcal{L}_v\beta'=\mathcal{L}_uf^*\beta'$,
$f^*i_{v'}d\beta=i_{u'}df^*\beta$, and $f^*(i_{v\wedge
v'}\phi)=i_{u\wedge u'}f^*\phi$. Using these relations and the
integrability of $L$, it follows that
$\Cour{(v,\beta),(v',\beta')}|_{f(Q)}\in \mathfrak{F}_f(L)$. The
map $\mathfrak{B}_f$ can be treated analogously.

Finally, for $L\in \mathrm{Dir}(M)_{f(Q)}$, one can check that
$\mathfrak{F}_f\circ \mathfrak{B}_f(L)= L$ if and only if
$\mathrm{pr}_{TM}(L)\subseteq Tf(Q)$, which holds since $L$ is
supported on $f(Q)$. Similarly, for $L\in \mathrm{Dir}(Q)$,
$\mathfrak{B}_f\circ \mathfrak{F}_f(L)=L$ if and only if
$\mathrm{ker}(df)\subseteq TQ\cap L$, which holds since $f$ is an
embedding.
\end{pf}

\subsection{Morphisms of Manin pairs}

Let $(E_i,A_i)$ be a Manin pair over $S_i$, $i=1,2$. We use the
pairing in $E_i$ to identify $E_i/A_i\cong A_i^*$, and denote by
$p_i:E_i \to A_i^*$ the natural projection.

\begin{defn}\label{def:morphismMP}
A \textbf{morphism} of Manin pairs $(E_1,A_1) \to (E_2,A_2)$ is a
Dirac structure $\D$ in $E_1\times \overline{E}_2$ with support on
the graph of a smooth map $J:S_1\to S_2$ (i.e., a Courant
algebroid morphism $E_1\to E_2$), such that the image of $\D$
under the projection $(p_1,p_2):E_1\times E_2 \to A_1^*\times
A_2^*$ is the graph of a bundle map $A_1^* \to J^*A_2^*$.
\end{defn}

Composition of these morphisms is by definition composition of relations
(see also Rem.\ \ref{rem:super}, where the morphisms are interpreted as maps,
and the composition becomes just composition of maps).

The next proposition gives a more explicit characterization of
morphisms of Manin pairs.

\begin{prop}\label{prop:equivMP}
Let $(E_i,A_i)$ be  a Manin pair over $S_i$, $i=1,2$, and let
$J:S_1\to S_2$ be a smooth map. A Dirac structure $\D$ in
$E_1\times \overline{E}_2$ with support on $\mathrm{graph}(J)$ is
a morphism of Manin pairs $(E_1,A_1)\to (E_2,A_2)$ if and only if
\begin{itemize}
\item[$i)$] $\D\cap A_1=\{0\}$,

\item[$ii)$] $\D \cap (A_1 \oplus J^*E_2)$ projects onto $J^*A_2$
under the natural projection $E_1\oplus J^*E_2\to J^*E_2$.
\end{itemize}
In other words, $i)$ and $ii)$ say that the projection $E_1\oplus
J^*E_2 \to J^*E_2$ restricted to $\D \cap (A_1 \oplus J^*E_2)$ is
an isomorphism onto $J^* A_2$.
\end{prop}

\begin{pf}
Let us consider the projection $p=(p_1,p_2): E_1\times E_2 \to
A_1^*\times A_2^*$, and let $R=p(\D)\subset A_1^*\oplus J^*
A_2^*$. We must show that conditions $i)$ and $ii)$ are equivalent
to $R$ being the graph of a bundle map $A_1^* \to J^*A_2^*$, or,
equivalently, that $R$ projects isomorphically onto $A_1^*$.

The projection of an element $r=p(k)\in R$ on $A_1^*$ is zero if
and only if $k\in \D\cap (A_1\times E_2)$. In this case, $r=0$ if
and only if its projection on $A_2^*$ is zero, which is equivalent
to the projection of $k$ on $E_2$ lying in $A_2$. Hence $R$
projects injectively on $A_1^*$ if and only if $ii)$ holds.  If
the projection of $R$ on $A_1^*$ is also onto, then $i)$ must
hold: given $a\in \D\cap A_1$ and any $\xi \in A_1^*$, there is a
$k\in \D$ with $p_1(k)=\xi$, and since $\D$ is isotropic, one has
$\SP{k,a}=\xi(a)=0$. This implies that $a=0$. Conversely, a
dimension count shows that the projection of $R$ on $A_1^*$ is
onto if and only if $\mathrm{rank}(K\cap(A_1\times
E_2))=\mathrm{rank}(A_2)$, which follows from $i)$ and $ii)$.
\end{pf}

\begin{numex} Let $S_1, S_2$ be manifolds and $f:S_1 \to
S_2$ a smooth map. Fix a closed 3-form $\phi\in \Omega^3(S_2)$,
and consider the Courant algebroids $\TS_1$, $\TS_2$, with
brackets defined by the 3-forms $f^*\phi$, $\phi$, respectively
(see Example \ref{ex:standard}). There is a natural Dirac
structure $K$ on $\TS_1\times \overline{\TS}_2$ supported on
$\mathrm{graph}(f)$ (c.f. \eqref{gds-basic}),
\begin{equation}\label{eq:Kstrong}
K_{(x,f(x))}=\{ ((X,f^*\beta),(df(X),\beta))\,|\, X\in T_xS_1,\;
\beta \in T^*_{f(x)}S_2 \}.
\end{equation}
Let $L_i \subset \TS_i$ be Dirac structures, so that $(\TS_i,L_i)$
is a Manin pair, $i=1,2$. Then property $i)$ in
Prop.~\ref{prop:equivMP} amounts to the condition
$\mathrm{ker}(df)\cap(TS_1\cap L_1)=\{0\}$, whereas $ii)$ says
that $f$ is a \textit{forward Dirac map} \cite{BR}. Hence $K$
defines a morphism of Manin pairs $(\TS_1,L_1)\to (\TS_2,L_2)$ if
and only if $f$ is a \textit{strong Dirac map} as in
\cite{ABM,BC,BC2} (see also Section \ref{subsec:twoclasses}). More
generally, one can fix a 2-form $\omega\in \Omega^2(S_1)$, define
the Courant bracket on $\TS_1$ using $f^*\phi + d\omega$, and
check that $K_{(x,f(x))}=\{
((X,f^*\beta-i_X\omega),(df(X),\beta))\,|\, X\in T_xS_1,\; \beta
\in T^*_{f(x)}S_2 \}$ is still a Dirac structure. Then $K$ defines
a morphism of Manin pairs if and only if $(f,\omega)$ is a strong
Dirac morphism as in \cite[Sec.~2.2]{ABM}.
\end{numex}

\begin{numrmk}[Super-geometric viewpoint]\label{rem:super}
The notion of morphism of Manin pairs in Def.~\ref{def:morphismMP}
has a natural interpretation in the framework of Gerstenhaber algebras, or equivalently
 in terms of odd Poisson structures,
as we now briefly outline.

As observed in \cite{Se2}, Manin pairs $(E,A)$ are in one-to-one
correspondence with principal $\mathbb{R}[2]$-bundles $P\to
A^*[1]$ (in the category of graded manifolds) equipped with a
$\mathbb{R}[2]$-invariant Poisson structure on $P$ of degree $-1$.
A morphism between $P_1\to A_1^*[1]$ and $P_2\to A_2^*[1]$ is
clear in this context: it is an $\mathbb{R}[2]$-equivariant
Poisson map. When such morphisms are expressed in terms of the
associated Manin pairs, one recovers
Definition~\ref{def:morphismMP}.

To see how this correspondence goes, given $P\to A^*[1]$, choose a
trivialization of $P$ to $A^*[1]\times \mathbb{R}[2]$, so that the
algebra of functions on $P$ becomes $\Gamma(\wedge A)[t]$, where
$t$ (the coordinate on $\mathbb{R}[2]$) has degree 2. The Poisson
structure on $P$ amounts to a Gerstenhaber bracket on
$\Gamma(\wedge A)[t]$. By invariance, it descends to a
Gerstenhaber bracket on $\Gamma(\wedge A)$, which defines a Lie
algebroid structure on $A$. Fixing this Lie algebroid structure,
the Gerstenhaber bracket on $\Gamma(\wedge A)[t]$ is determined by
$[t,t]\in \Gamma(\wedge^3 A)$ and $[t,a]\in \Gamma(\wedge A)$, for
all $a\in \wedge A$. By setting $\chi:=[t,t]$ and
$d_{A^*}:=[t,\cdot]:\Gamma(\wedge^\bullet A)\to
\Gamma(\wedge^{\bullet -1}A)$, we define a Lie quasi-bialgebroid.
We hence set $E=A\oplus A^*$ to be the double Courant algebroid,
so that $(E,A)$ is a Manin pair. Different trivializations of $P$
amount to different choices of isotropic complements of $A$ in
$E$, so this procedure establishes the desired correspondence.

This correspondence can also be viewed more instrinsically,
without the choice of splittings or trivializations.

First recall from \cite{Se} (concluding work of Vaintrob,
Roytenberg \cite{Roy0} and Weinstein) that Courant algebroids $E$
are equivalently described as degree 2 symplectic non-negatively
graded manifolds $\mathcal{M}$ equipped with a function $\Theta$
of degree 3 satisfying $\{\Theta,\Theta\}=0$ (c.f. \cite{Roy2}).
Dirac structures in $E$ with support on a submanifold correspond
to the Lagrangian submanifolds of $\mathcal{M}$ on which the
degree 3 function $\Theta$ vanishes. This motivates the definition
of morphism of Courant algebroids in
Section~\ref{subsec:diracsupp}.

The Manin
pair associated with $P\to A^*[1]$ is obtained as follows. The
Poisson structure on $P$ is a function $H$ on $T^*[2]P$, quadratic
on the fibres, of total degree 3 and such that $\{H,H\}=0$. This
function descends to the symplectic reduction
$Y=T^*[2]P\slash\hspace{-1mm}\slash_1\mathbb{R}[2]$ at moment
value $1$, hence $Y$ determines a Courant algebroid $E$ (see
\cite{BCMZ} for more on the reduction of Courant algebroids from
this perspective). The map $Y\to A^*[1]$ (which is a Lagrangian
fibration) gives rise to a map $E\to A^*$, whose kernel determines
$A\subset E$, defining a Manin pair $(E,A)$. Given $P_1\to
A_1^*[1]$ and $P_2\to A_2^*[1]$ and a morphism $\psi: P_1\to P_2$,
since $\psi$ is a Poisson map, its graph $\Gamma_\psi$ is
coisotropic in $P_1 \times \overline{P}_2$ (where $\overline{P}_2$
is $P_2$ with the opposite Poisson structure). Then the conormal
bundle $N^*[2]\Gamma_\psi$ is Lagrangian in $T^*[2](P_1\times
P_2)$, and the function $H$ coming from the Poisson structure on
$P_1\times \overline{P}_2$ vanishes on it. The reduction of
$N^*[2]\Gamma_\psi$ to the symplectic quotient $Y_1\times
\overline{Y}_2$ defines a Lagrangian submanifold, which
corresponds to the Dirac structure $\D$ of
Definition~\ref{def:morphismMP}.
\end{numrmk}

\section{Morphisms of Manin pairs and Hamiltonian spaces}\label{sec:ham}

\subsection{Hamiltonian spaces}\label{subsec:ham}

Let us consider a Manin pair $(E, A)$
over a manifold $S$. The following
is the main definition of this paper.

\begin{defn}\label{Hamiltonian-space}
A \textbf{Hamiltonian space} for $(E, A)$ is a manifold $X$
together with a morphism of Manin pairs
$(\mathbb{T}X,TX)\to(E,A)$, where $\mathbb{T}X$ is the standard
Courant algebroid of $X$.
\end{defn}
 Hamiltonian spaces are denoted
by triples $(X,J,\D)$  (where the map $J:X\to S$ and the Dirac
structure $\D$ are as in Definition \ref{def:morphismMP}), and $J$
is the \textbf{moment map}.

A \textbf{morphism} between two Hamiltonian spaces $(X,J,\D)$ and
$(X',J',\D')$ is a  smooth map $f:X \to X'$ such that
$J(x)=J'(f(x))$, $\forall x\in X$, and $\Gamma_f\circ \D=\D'$,
where $\Gamma_f$ is defined in \eqref{gds-basic}, Example
\ref{ex:gammaf}, and $\circ$ denotes \textit{composition of
relations}, i.e.,
\begin{eqnarray*}
(\D')_{(f(x),J(x))}&=&\{((u',\alpha'),e)\;|\; (u',\alpha')\in
(\TX')_{f(x)},\, e\in E_{J(x)}, \\
&& \mbox{ and } \exists u\in (TX)_x \mbox{ s.t. } u'=df(u),\,
((u,f^*(\alpha')),e)\in (\D)_{(x,J(x))} \},\; \forall x\in X,
\end{eqnarray*}
see e.g. \cite{BR,IX}. (Note that $(f(x),J(x))\in
\mathrm{graph}(J')$.) The \textbf{category of Hamiltonian spaces}
associated with the Manin pair $(E,A)$ is denoted by $\M(E,A)$.

Notice that the projection
$\mathbb{T}X\times E \to E$ restricts to an isomorphism
\begin{equation}\label{eq:iso}
\D_{(x,J(x))}\cap ((TX\oplus \{ 0\})_x\times
E_{J(x)})\stackrel{\sim}{\longrightarrow}A_{J(x)},\;\;\;\; \forall
x\in X.
\end{equation}



\begin{numex}\label{ex:canonical}
We present two examples of Hamiltonian spaces canonically
associated with any pair $(E,A)$ over $S$. Let $\Delta \subset
S\times S$ be the diagonal, and consider the subbundle $\D\subset
(\TS\times E)|_\Delta$ given by
$$
\D_{(x,x)}:=\{ ((\rho(a),-\beta), a + \rho^*\beta)\;|\; a\in A_x,
\, \beta \in (T^*S)_x \} \subset (\TS)_x\times E_x,\;\; x\in S.
$$
This is a Dirac structure in $\TS\times E$ supported on $\Delta$,
and it defines a Hamiltonian space $(X,J,\D)$ for $(E,A)$, where
$X=S$ and $J=\mathrm{Id}$. Indeed, a direct check using
\eqref{eq:rhorho*} shows that $\D$ is isotropic, and a dimension
count gives that it has maximal rank. Condition $d1)$ in
Def.~\ref{def:dirac} is a consequence of \eqref{eq:rhorho*}. To
verify the integrability condition $d2)$, it suffices to consider
sections of $\TS\times E$ of the form
$\widetilde{a}=((\rho(a),0),a)$, $a\in \Gamma(A)$, and
$\widetilde{\beta}=((0,-\beta),\rho^*\beta)$, $\beta \in
\Omega^1(S)$, and check that
$$
\Cour{\widetilde{a},\widetilde{a}'}|_\Delta,\;\;
\Cour{\widetilde{\beta},\widetilde{\beta}'}|_\Delta,\;\;
\Cour{\widetilde{a},\widetilde{\beta}}|_\Delta \in \Gamma(\D),
$$
see Remark~\ref{rmk:bases}. The first case follows directly from
$c4)$. For the second case, \eqref{eq:prop2} gives
$\Cour{\widetilde{\beta},\widetilde{\beta}'}=0$. For the last
case, note that
$\Cour{\widetilde{a},\widetilde{\beta}}=((0,-\mathcal{L}_{\rho(a)}\beta),\Cour{a,\rho^*\beta})$.
It is immediate to check (using $c4)$ and \eqref{eq:rhorho*}) that
$\langle
\Cour{\widetilde{a},\widetilde{\beta}},\widetilde{\beta}'\rangle=0$,
$\forall \beta'$. On the other hand, using $c3)$, for all
$\widetilde{a}'$ we have
$$
\langle \Cour{\widetilde{a},\widetilde{\beta}},\widetilde{a}'
\rangle = -i_{\rho(a')}\mathcal{L}_{\rho(a)}\beta + \langle a',
\Cour{a,\rho^*\beta} \rangle =
-i_{\rho(a')}\mathcal{L}_{\rho(a)}\beta -\langle
\Cour{a,a'},\rho^*\beta \rangle +
\mathcal{L}_{\rho(a)}i_{\rho(a')}\beta = 0,
$$
since $\langle \Cour{a,a'},\rho^*\beta \rangle =
i_{[\rho(a),\rho(a')]}\beta =
\mathcal{L}_{\rho(a)}i_{\rho(a')}\beta -
i_{\rho(a')}\mathcal{L}_{\rho(a)}\beta$. Hence $\langle
\Cour{\widetilde{a},\widetilde{\beta}}|_\Delta,
\Gamma(K)\rangle=0$, and since $\D$ is maximal isotropic,
$\Cour{\widetilde{a},\widetilde{\beta}}|_\Delta\in \Gamma(K)$.

Similarly, let $\iota: \mathcal{O}\hookrightarrow S$ be an orbit
of $A$ (i.e., an integral leaf of $\rho(A)\subseteq TS$), and
$\D\subset (\T\mathcal{O}\times E)|_{\mathrm{graph}(\iota)}$ be
given by $ \D_{(x,\iota(x))}=\{ ((\rho(a),-\iota^*\beta), a +
\rho^*\beta)\;|\; a\in A_x, \,\beta \in (T^*S)_x \}.$ Then
$(X,J,\D)$, where $X=\mathcal{O}$ and $J=\iota$, is a Hamiltonian
space for $(E,A)$.
\end{numex}

The definition of Hamiltonian space implies the following
properties: 

\begin{prop}\label{induced-action}
Let $(X,J,\D)$ be a Hamiltonian space for $(E, A)$. Then:
\begin{itemize}
\item[{\it i)}] For any $a\in A_{J(x)}$, there exists a unique
$u\in T_xX$ such that $((u,0),a)\in \D_{(x,J(x))}$. Moreover, this
induces a Lie algebroid action $\rho_X:J^*A\to TX$ of $A$ on
$J:X\to S$ (i.e., $\rho_X$ is a smooth bundle map so that
$dJ(\rho_X(a)_x)=\rho(a_{J(x)})$ and the induced map
$\rho_X:\Gamma(A)\to \mathfrak{X}(X)$ preserves Lie brackets.)

\item[{\it ii)}] For any $\alpha \in T_x^*X$, there exists $e\in
E_{J(x)}$ and $u\in T_xX$ such that $((u,\alpha),e)\in
\D_{(x,J(x))}$. In addition, $\rho_X^*(\alpha_x)=0$ if and only if
$e\in A_{J(x)}$, for $x\in X$.

\item[{\it iii)}] Let $A'$ be any complement of $A$ in $E$, so
that $E=A\oplus A'$. Let us identify $A'$ with
$\{((0,0),(0,a'))\;|\; a'\in A' \}$ in $\TX\times E$. Then $\D\cap
A'=\{0\}$.


\end{itemize}

\end{prop}

\begin{pf}
The first assertion in $i)$ follows from $i)$ and $ii)$ in Prop.
\ref{prop:equivMP}. Note that
$\rho_X:J^*A\to TX$ is defined by the inverse of \eqref{eq:iso}
followed by projection on $TX$, so it is smooth, and it preserves
brackets as a result of the integrability of $\D$. The property
$dJ(\rho_X(a))=\rho(a)$, $a\in A_{J(x)}$, is a direct consequence
of $d1)$ in Def.~\ref{def:dirac}.

Consider the projection $p : \D\to T^*X$. Using $ii)$ in Prop.
\ref{prop:equivMP}, one can deduce that $\mbox{ker}(p)
=\D\cap (E\times (TX\oplus \{ 0\}))\cong A$. A dimension count
shows that $p$ is surjective, so the first statement in $ii)$
follows. For the second assertion, we consider the pairing of
$((u,\alpha),e) \in \D$ with elements of the form
$((\rho_X(a),0),a)$, for $a\in J^*A$ (which are necessarily in
$\D$ by definition of $\rho_X$), and use that $A$ is maximal
isotropic.


Property $iii)$ is a direct consequence of  $ii)$ in Prop.\ \ref{prop:equivMP} and $A\cap
A'=\{0\}$.
\end{pf}

\subsection{Two different characterizations}\label{subsec:twoclasses}

We now see how the intrinsic notion of Hamiltonian space for a
Manin pair $(E,A)$, discussed in Section~\ref{subsec:ham}, has
more familiar incarnations once extra (noncanonical) choices are
made.

\subsubsection*{Quasi-Poisson geometry}

Let $(E,A)$ be a Manin pair over $S$.  The first possible
noncanonical choice we consider is that of an isotropic splitting
$j$ as in \eqref{eq:j}, defining a Lie quasi-bialgebroid structure
on $A$ with cobracket $F_j:\Gamma(A)\to \wedge^2\Gamma(A)$, bundle
map $\rho_{A^*}^j:A^*\to TM$ and  3-tensor $\chi_j\in
\Gamma(\wedge^3A)$ defined as in \eqref{eq:liequasi}.

Let $(X,J,\D)$ be a Hamiltonian space for the Manin pair $(E,A)$.
One can adapt the construction in \cite{BCS} to show that, once
$j$ is fixed, one naturally obtains a bivector field $\Pi_X^j \in
\mathfrak{X}^2(X)$ as follows. Conditions $i)$ and $ii)$ in Prop.\ \ref{prop:equivMP} imply that, for all $x\in X$,
$\D_{(x,J(x))} \cap ((TX)_x\times (A^*)_{J(x)})=\{0\}$, i.e.,
$$
{\widehat{J}}^*\D\cap (TX \oplus J^*A^*)=\{0\},
$$
where $\widehat{J}=\mathrm{Id}\times J: X\to
\mathrm{graph}(J)\subset X\times S$. Hence $\widehat{J}^*\D$ is
the graph of a skew-symmetric bundle map $T^*X \oplus J^*A \to TX
\oplus J^*A^*$, determined by an element $\Lambda^j \in
\Gamma(\wedge^2(TX\oplus J^*A^*))$. The bivector field $\Pi_X^j$
is the component of $\Lambda^j$ in $\Gamma(\wedge^2TX)$: given
$\alpha\in (T^*X)_x$, $i_\alpha\Pi_X^j$ is the only element in
$(TX)_x$ satisfying
\begin{equation}\label{eq:PiX}
((i_\alpha\Pi_X^j,\alpha),(0,-\rho_X^*\alpha)) \in \D_{(x,J(x))},
\end{equation}
where $\rho_X:J^*A\to TX$ is the action map given in
Prop.~\ref{induced-action}, part $i)$, see \cite{BCS}. We have the
following alternative characterization of $\Pi_X^j$:
\begin{lem}\label{lem:piX}
Viewing $A^*\subset E$ via $j$, we have that
$\mathrm{graph}(\Pi_X^j)=\D\circ A^*$, where $\circ$ denotes
composition of relations.
\end{lem}
\begin{pf}
The composition of $\D_{(x,J(x))}$ and $(A^*)_{J(x)}$ gives
\[
(\D \circ A^* )_x=\{ (u,\alpha )\in (\TX)_x \, | \, \exists\, e\in
A^*_{J(x)}\, \mathrm{ with } \; ((u,\alpha ),e)\in
\D_{(x,J(x))}\}.
\]
This is a Lagrangian subspace of $(TX\oplus T^*X)_x$, so it
suffices to show that $\mathrm{graph}(\Pi_X^j)\subseteq \D\circ
A^*$ pointwise. But this is a direct consequence of
\eqref{eq:PiX}.
\end{pf}

For an arbitrary Lie quasi-bialgebroid $(A, F ,\chi, \rho_{A^*})$
over $S$, there is a notion of {\bf Hamiltonian space} (see e.g.
\cite{ILX}) defined by quadruples $(X,\Pi_X,J,\rho_X)$, where $X$
is a manifold, $\Pi_X\in \mathfrak{X}^2(X)$, $J:X \to S$ is a
smooth map (the \textit{moment map}), and $\rho_X:J^*A\to TX$ is a
Lie algebroid action, satisfying the following compatibility
conditions:
\begin{align}
& \frac{1}{2}[\Pi_X,\Pi_X]=\rho_X(\chi), \label{eq:qpois1}\\
& \mathcal{L}_{\rho_X(a)}\Pi_X=\rho_X(d_{*}a),\; a\in \Gamma(A),\label{eq:qpois2}\\
& \Pi_X^\sharp J^*=\rho_X \rho_{A^*}^*,\label{eq:qpois3}
\end{align}
where $d_{*}$ is the quasi-differential on $\Gamma(\wedge A)$
determined by $F$ and $\rho_{A^*}$. These Hamiltonian spaces are
also referred to as \textit{Hamiltonian quasi-Poisson spaces}, and
they form a category in which morphisms between
$(X,\Pi_X,J,\rho_X)$ and $(X',\Pi_{X'},J',\rho_X')$ are smooth
equivariant maps $f:X\to X'$ such that $f_*\Pi_X=\Pi_{X'}$ and
$J'\circ f=J$.

Given a Manin pair $(E,A)$ together with a splitting $j$, we
denote by $\M_j(E,A)$ the category of Hamiltonian spaces
associated with the Lie quasi-bialgebroid
$(A,F_j,\chi_j,\rho_{A^*}^j)$ defined by $j$. The following result
generalizes \cite[Thm.~4.1]{BCS} and is proven similarly:

\begin{prop}\label{prop:quasi} Given a Hamiltonian space $(X,J,\D)$ for $(E,A)$ and
a splitting $j$, then $(X,\Pi_X^j,J,\rho_X)$ is a Hamiltonian
space for the Lie quasi-bialgebroid determined by $j$; conversely,
given a Hamiltonian space $(X,\Pi_X,J,\rho_X)$ for a Lie
quasi-bialgebroid $(A,F,\chi,\rho_{A^*})$, the bundle
\begin{equation}\label{eq:H}
\D
:=\{((\rho_X(a)+i_\alpha\Pi_X,\alpha),(a,-\rho_X^*(\alpha)))\;|\;
a\in J^*A, \, \alpha \in T^*X\}
\end{equation}
over $\mathrm{graph}(J)$ makes $(X,J,\D)$ into a Hamiltonian space
for $(E,A)$, where $E=A\oplus A^*$ is the double Courant
algebroid. Moreover, these constructions define an isomorphism of
Hamiltonian categories $\M(E,A)\cong \M_j(E,A)$.
\end{prop}

\begin{numex}[Quasi-Poisson $D/G$-valued moment maps]
If $(E=\mathfrak{d}_S, A=\frakg_S)$ is the Manin pair over $S=D/G$
of Example~\ref{ex:dressing}, and if $j$ is a splitting of the
Manin pair $(\frakd,\frakg)$ (i.e.,
$\mathfrak{h}=j(\mathfrak{g}^*)$ is an isotropic complement of
$\frakg$ in $\frakd$), then the Hamiltonian spaces in $\M_j(E,A)$
are exactly the Hamiltonian quasi-Poisson spaces with $D/G$-valued
moment maps considered in \cite{AK} (c.f. \cite[Sec.~5]{BC2}). For
the canonical Hamiltonian space $(S,\mathrm{Id},\D)$ of
Example~\ref{ex:canonical}, the associated bivector field given by
Prop.~\ref{prop:quasi} agrees with the one in \cite[Sec.~3.5]{AK}
(c.f. \cite[Sec.~5.1]{BC2}); it is the bivector field determined
by the Lie quasi-bialgebroid $\frakg_S, \mathfrak{h}_S \subset
\frakd_S$ over $S$ (c.f. \cite{ILX}, \cite[App.~4]{BC2}) .
\end{numex}

\subsubsection*{Dirac geometry}

To define a second possible realization of the Hamiltonian
category associated with a Manin pair $(E,A)$, we assume that $E$
is an \textit{exact Courant algebroid}, i.e., that the sequence
\begin{equation}\label{eq:exact}
0\to T^*S \stackrel{\rho^*}{\longrightarrow} E
\stackrel{\rho}{\longrightarrow} TS \to 0
\end{equation}
is exact. One can always choose (in a noncanonical way) an
isotropic splitting $s:TS\to E$, which induces an identification
of $E$ with $\TS=TS\oplus T^*S$, where the latter is equipped with
the Courant algebroid structure of Example~\ref{ex:standard}, with
closed 3-form $\phi_S^s\in \Omega^3(S)$ defined by
\begin{equation}\label{eq:phiS}
\phi_S^s(v_1,v_2,v_3):= \SP{s(v_1),\Cour{s(v_2),s(v_3)}}.
\end{equation}
Under this identification, $A\subset E$ defines a Dirac structure
$L^s_S$ on $S$ (with respect to $\phi_S^s$),
\begin{equation}\label{eq:LS}
L_S^s=\{ (\rho(a),s^*(a))\,|\, a\in A\} \subset \TS.
\end{equation}

We recall that there is a general notion of \textbf{Hamiltonian
space} associated with any Dirac manifold $(S,L_S,\phi_S)$
\cite{BC,BC2}, where we view $S$ as the target of the moment map.
These Hamiltonian spaces are given by smooth maps $J:X\to S$,
where $X$ is a manifold equipped with a Dirac structure $L_X$ for
which $J$ is a \textbf{strong Dirac map} \cite{ABM,BC,BC2}; that
is, $J:X\to S$ is a forward Dirac map \cite{BR}, $L_X$ is
integrable with respect to $J^*\phi_S$, and the following
transversality condition holds: $\mathrm{ker}(dJ)\cap (L_X\cap
TX)=\{0\}$. We denote such Hamiltonian spaces by triples
$(X,L_X,J)$. A morphism between Hamiltonian spaces $(X,L_X,J)$ and
$(X',L_{X'},J')$ is a forward Dirac map $f:X\to X'$ such that
$J'\circ f=J$, and the category of Hamiltonian spaces associated
with $S$ is denoted  by $\M(S,L_S,\phi_S)$.

Given a Manin pair $(E,A)$ over a manifold $S$ such that $E$ is
exact, and letting $s:TS \to E$ be an isotropic splitting of
\eqref{eq:exact}, we have the associated category of Hamiltonian
spaces $\M_s(E,A):=\M(S,L_S^s,\phi_S^s)$. Let
$\widehat{J}:=(\mathrm{Id},J):X\to X \times S$, and consider the
maps $\mathfrak{F}_{\widehat{J}}$ and $\mathfrak{B}_{\widehat{J}}$
from Prop.~\ref{prop:correspondence} (with $Q=X$, and $M=X\times
S$ equipped with 3-form $0\times \phi_S^s$).

\begin{them}\label{thm:dirac}
Given a Hamiltonian space $(X,L_X,J)$ for $(S,L_S^s,\phi_S^s)$, then
$(X,J,\D=\mathfrak{F}_{\widehat{J}}(L_X))$ is a Hamiltonian space for
$(E,A)$ and, given $(X,J,\D)$ a Hamiltonian space for
$(E,A)$, then $(X,L_X=\mathfrak{B}_{\widehat{J}}(\D),J)$ is a
Hamiltonian space for $(S,L_S^s,\phi_S^s)$. Moreover, these procedures
establish an isomorphism of categories
$\M_s(E,A)\cong \M(E,A)$.
\end{them}

\begin{pf}
Given $(X,L_X,J)$, we consider
\begin{equation}\label{eq:FLX}
\D:= \mathfrak{F}_{\widehat{J}}(L_X)= \{ ((u,\alpha
-J^*\beta),(dJ(u), \beta ))\, | \, (u,\alpha )\in L_X,\, \beta \in
T^*S \},
\end{equation}
which is a Dirac structure in $\mathbb{T}X\times \mathbb{T}S$
supported on $\mathrm{graph}(J)$ (integrable with respect to
$0\times \phi_S^s$). In order to show that $\D$ defines a
Hamiltonian space for the Manin pair $(E,A)=(\mathbb{T} S,L_S^s)$,
we must show that it satisfies conditions $i)$ and $ii)$ in Prop.\ \ref{prop:equivMP}.

To prove $i)$, note that if $((u,0),(0,0))\in \D$, then $(u,0)\in
L_X\cap TX$ and $u\in \mbox{ker}(dJ)$. Hence $u=0$ since
$\mathrm{ker}(dJ)\cap (L_X\cap TX)= \{ 0 \}$. To prove $ii)$, let
$((u,0),(dJ(u),\beta ))\in \D$. Then $(u,J^*\beta )\in L_X$ and,
since $J$ is forward Dirac map, we have $(dJ(u),\beta )\in L_S$.
On the other hand, if $(v,\beta)\in L_S$ then, using that
$J:(X,L_X)\to (S, L_S)$ is a forward Dirac map, there exists $u\in
TX$ such that $dJ(u)=v$ and $(u,J^*\beta )\in L_X$. Thus,
$((u,0),(dJ(u),\beta ))\in \D$.

Let now $(X,J,\D)$ be a Hamiltonian space for the Manin pair
$(E,A)=(\mathbb{T}S,L_S^s)$. Then $L_X:=
\mathfrak{B}_{\widehat{J}}(\D)$ is a Dirac structure on $X$
(integrable with respect to $J^*\phi _S^s$). Explicitly,
\begin{equation}\label{Induced-Dirac}
L_X=\{ (u,\alpha +J^*\beta )\, | \, ((u,\alpha ),(dJ(u), \beta
))\in \D \}.
\end{equation}
To see that $J$ is a forward Dirac map, it suffices to show that
$L_S\subseteq \mathfrak{F}_J(L_X)$ at each point $J(x)$. So take
$(v,\beta)\in L_S$ at $J(x)$. By $ii)$, there exists $u \in
(TX)_x$ such that $((u,0),(v,\beta))\in \D$, with $v=dJ(u)$. It
follows from \eqref{Induced-Dirac} that $(u,J^*\beta)\in L_X$, so
$(v,\beta)\in \mathfrak{F}_J(L_X)$. It remains to check that
$\mathrm{ker}(dJ)\cap (L_X\cap TX )=\{0\}$. From
\eqref{Induced-Dirac}, we see that $u\in \mathrm{ker}(dJ)\cap
(L_X\cap TX)$ if and only if there exists $\beta\in T^*S$ with
$((u,-J^*\beta),(0,\beta))\in \D$. Now notice that any element of
the form $((0,-J^*\beta),(0,\beta))$ is in $\D$ (just check that
the pairing of an element of this form with any element in $\D$
must vanish using that $\D$ is supported on the graph of $J$). It
follows that $((u,0),(0,0))\in \D$, hence $u=0$ by $i)$.

Finally, note that $f:(X,L_X)\to (X',L_{X'})$ is a forward Dirac
map if and only if $L_{X'}=\Gamma_f \circ L_X$, and, in case
$J'\circ f=J$, this is equivalent to $\D'=\Gamma_f\circ \D$. So
morphisms in $\M(E,A)$ and $\M_s(E,A)$ are naturally identified.
\end{pf}

\begin{numex}[Dirac geometric $D/G$-valued moment maps]
When  $(E=\mathfrak{d}_S, A=\frakg_S)$ is the Manin pair over
$S=D/G$ of Example~\ref{ex:dressing}, and $s$ is an isotropic
splitting of \eqref{eq:exact}, then the Hamiltonian spaces in
$\M_s(E,A)$ coincide with the ones considered in
\cite[Sec.~4]{BC2}, and particular examples include the $G$-valued
moment maps of \cite{AMM} (c.f. \cite{ABM,BC,BC2}). For the
canonical Hamiltonian space $(S,\mathrm{Id},\D)$ of
Example~\ref{ex:canonical}, the associated Dirac structure on $S$
induced by $s$ is simply $L_S^s$ given in \eqref{eq:LS}; the Dirac
structure on a dressing $\mathfrak{g}$-orbit $\mathcal{O}$ is its
presymplectic structure as a leaf of the Dirac manifold
$(S,L^s_S)$.
\end{numex}

\subsection{Equivalence}

Let $(E,A)$ be a Manin pair over a manifold $S$, and assume that
$E$ is an exact Courant algebroid. Let us fix both types of extra
choices considered in Section~\ref{subsec:twoclasses}: an
isotropic splitting $j:A^*\to E$ (making $A$ into a Lie
quasi-algebroid $(A,F,\chi,\rho_{A^*})$ given by
\eqref{eq:liequasi}) as well as an identification $E\cong TS\oplus
T^*S$, induced by $s:TS\to E$, where the Courant algebroid
structure on $\TS$ is with respect to the 3-form $\phi_S$ given in
\eqref{eq:phiS}. (We are simplifying the notations by omitting the
dependence on $j$ and $s$.) Under this identification, $A\subset E$
defines a Dirac structure $L_S$ on $S$ \eqref{eq:LS}, whereas the subbundle
$A^*\subset E$ (via $j$) defines a transverse almost Dirac structure $C_S$.
As discussed in
Section~\ref{subsec:twoclasses}, each choice leads to a category
of Hamiltonian spaces, denoted by $\M_j(E,A)$ and $\M_s(E,A)$.
Since Prop.~\ref{prop:quasi} and Thm.~\ref{thm:dirac} say that
both categories are different realizations of the same category
$\M(E,A)$, we immediately obtain:


\begin{them}\label{thm:equiv}
There is an isomorphism of categories $\M_j(E,A)\cong \M_s(E,A)$
as follows: given a Hamiltonian space $(X,\Pi_X,J,\rho_X)$ for the
Lie quasi-bialgebroid $(A,F,\chi,\rho_{A^*})$, then the triple
$(X,L_X,J)$ is a Hamiltonian space for the Dirac manifold
$(S,L_S,\phi_S)$, where
\begin{equation}\label{eq:LX}
L_X:=\{ (\rho_X(a) + i_\alpha \Pi _X, J^*s^*(a) + (\mathrm{Id}
- J^*\overline{\rho}^*\rho_X^*)\alpha)\,|\, a\in A,\, \alpha \in
T^*X \},
\end{equation}
with $\overline{\rho}=j^*s:TS\to A$. Conversely, given a
Hamiltonian space $(X,L_X,J)$ for $(S,L_S,\phi_S)$, then the
composition of relations $\mathfrak{F}_{\widehat{J}}(L_X)\circ
A^*\subset \TX$ defines a bivector field $\Pi_X\in \mathfrak{X}^2(X)$,
\begin{equation}\label{eq:piX2}
\mathrm{graph}(\Pi_X)=\mathfrak{F}_{\widehat{J}}(L_X)\circ A^*,
\end{equation}
making $(X,\Pi_X,J,\rho_X)$ into a Hamiltonian quasi-Poisson space
for $(A,F,\chi,\rho_{A^*})$.
\end{them}

\begin{pf}
The proof follows by combining the constructions identifying both
categories $\M_j(E,A)$ and $\M_s(E,A)$ with $\M(E,A)$: the
expression for the Dirac structure \eqref{eq:LX} is obtained
directly as $\mathfrak{B}_{\widehat{J}}(\D)$, where $\D$ is given
by \eqref{eq:H}. Conversely, the expression for the quasi-Poisson
bivector field \eqref{eq:piX2} follows from Lemma \ref{lem:piX}.
\end{pf}

Let us denote by $\mathcal{M}(E,A)$ the subcategory of $\M(E,A)$
consisting of Hamiltonian spaces satisfying the additional
condition that the natural projection of $\D$ on $TX$ is onto.
Then the isomorphism in Thm.~\ref{thm:equiv} restricts to an
isomorphism of subcategories,
$$
\mathcal{M}_j(E,A)\cong \mathcal{M}(E,A) \cong \mathcal{M}_s(E,A),
$$
where Hamiltonian spaces $(X,\Pi_X,J,\rho_X)$ in $\mathcal{M}_j(E,A)$ satisfy the extra condition
$TX=\{\rho_X(a)+i_\alpha\Pi_X\;|\; a\in J^*A,\,
\alpha \in T^*X  \}$, and
Hamiltonian spaces $(X,L_X,J)$ in $\mathcal{M}_s(E,A)$ are such that $L_X=\mathrm{graph}(\omega_X)$
for a given 2-form $\omega_X\in \Omega^2(X)$.

\begin{numrmk} The bivector field $\Pi_X$ \eqref{eq:piX2} admits another description, following \cite{ABM,BC2}.
Using \eqref{eq:FLX}, the expression in \eqref{eq:piX2} can be written explicitly as
$$
\mathrm{graph}(\Pi_X)_x=\{ (u,\alpha-J^*\beta)\in (\TX)_x \;|\;
(dJ(u),\beta)\in C_S,\; (u,\alpha) \in L_X \}.
$$
On the other hand, as shown in \cite{ABM}, the pull-back image (in
the Dirac geometric sense \cite{BR}) of $C_S$ under $J$ defines an
almost Dirac structure $C_X$ on $X$ transverse to $L_X$. The
splitting $\TX=L_X\oplus C_X$ defines a Lie quasi-bialgebroid, and
it naturally induces a bivector field $\Pi$ on $X$ (see e.g.
\cite{ILX} or \cite[App.~4]{BC2}). By \cite[Prop.~1.16]{ABM}, such
$\Pi$ is given by
$$
\mathrm{graph}(\Pi)_x=\{ (u,\alpha-J^*\beta)\in (\TX)_x \;|\;
(u,J^*\beta)\in C_X,\; (u,\alpha) \in L_X \}.
$$
Since $C_X$ is the pull-back of $C_S$, it immediately follows that
$\mathrm{graph}(\Pi_X)\subseteq \mathrm{graph}(\Pi)$. This implies that $\Pi_X=\Pi$.
\end{numrmk}

\begin{numex}
For the particular Manin pair  $(E=\mathfrak{d}_S, A=\frakg_S)$
over $S=D/G$ of Example~\ref{ex:dressing}, Thm.~\ref{thm:equiv}
recovers the equivalence proven in \cite[Sec.~6.3]{BC2}; this
result gives, as a special case, the equivalence between the two
formulations of $G$-valued moment maps in \cite{AMM} and
\cite{AKM} (see \cite[Sec.~10]{AKM}, \cite[Sec.~5.4]{ABM},
\cite[Sec.~3.5]{BC} for proofs).
\end{numex}

\subsection{Poisson algebras and moment map reduction}

Let $(X,J,\D)$ be a Hamiltonian space for a Manin pair $(E,A)$
over $S$. Following Prop.~\ref{induced-action}, part $i)$, let
$\rho_X$ denote the induced Lie algebroid action of $A$ on $J:X\to
S$.

A function $f\in C^\infty(X)$ is called \textbf{admissible} if
there exists a vector field $u_f\in \mathfrak{X}(X)$ satisfying
\begin{equation}\label{eq:hamvec}
((u_f,df)_x,0)\in \D_{(x,J(x))},\;\;\; \forall \; x\in X.
\end{equation}
By $i)$ in Prop.\ \ref{prop:equivMP}, $u_f$ is
\textit{uniquely} defined by condition \eqref{eq:hamvec}. Since
$\D$ is supported on graph($J$), $u_f$ satisfies $dJ(u_f)$=0. It
is simple to check that the set of admissible functions
$C^\infty(X)_{\mathrm{adm}}$ is a subalgebra of $C^\infty(X)$. We
define a bracket on $C^\infty(X)_{\mathrm{adm}}$ by
\begin{equation}
\{f,g\}:=\Lie_{u_f}g.
\end{equation}

\begin{lem}\label{lem:pois}
$(C^\infty(X)_{\mathrm{adm}},\{\cdot,\cdot\})$ is a Poisson
algebra.
\end{lem}

\begin{pf}
Let $f$ and $g$ be admissible functions. The fact that $\D$ is
isotropic implies that $df(u_g)=-dg(u_f)$, i.e., $\{\cdot,\cdot\}$
is skew-symmetric. Using the integrability of $\D$ with respect to
the product Courant bracket on $\TX\times E$ (condition $d2)$ in
Def.~\ref{def:dirac}), one can check that $\{f,g\}$ is admissible;
in fact, $u_{\{f,g\}}=[u_f,u_g]$. This last property also proves
the Jacobi identity for $\{\cdot,\cdot\}$.
\end{pf}

\begin{lem}\label{lem:invadm} A function $f$ satisfies $\Lie_{\rho_X(a)}f=0$
for all $a\in A_{J(x)}$ if and only if there exists $u_x\in
(TX)_x$ such that $((u_x,(df)_x),0)\in \D_{(x,J(x))}$.
\end{lem}

\begin{pf}
Since $((\rho_X(a),0),a)\in \D$, if we assume that
$((u_x,(df)_x)_,0)\in \D_{(x,J(x))}$ and use that $\D$ is
isotropic, it follows that $df(\rho_X(a))=\Lie_{\rho_X(a)}f=0$.
Conversely, suppose that $\Lie_{\rho_X(a)}f=0$, $a\in A_{J(x)}$.
From Prop.~\ref{induced-action}, part $ii)$, there exists $u'\in
T_xX$, $a'\in A_{J(x)}$ such that $((u',df),a')\in \D$. Since
$((\rho_X(a'),0),a')\in \D_{(x,J(x))}$, we have that
$((u_x,(df)_x),0)\in \D_{(x,J(x))}$, for $u_x=u'-\rho_X(a')$.
\end{pf}

A function $f\in C^\infty(X)$ is called \textbf{$A$-invariant} if
$\Lie_{\rho_X(a)}f=0$, $\forall a\in \Gamma(A)$, and the set of
all $A$-invariant functions is denoted by $C^\infty(X)^A$.

\begin{prop}\label{prop:poisalg}
A function $f\in C^\infty(X)$ is $A$-invariant if and only if it
is admissible, and therefore $C^\infty(X)^A$ is a Poisson algebra.
\end{prop}

\begin{pf}
From Lemma~\ref{lem:invadm}, we know that $f$ is $A$-invariant if
and only if, at each $x\in X$, there exists  $u_x\in (TX)_x$ such
that $((u_x,(df)_x),0)\in \D_{(x,J(x))}$. It remains to check that
$u_f:X\to TX$, $x\mapsto u_x$, is a smooth vector field. To see
that, fix any vector subbundle $A'\subset E$ such that $E=A\oplus
A'$, which defines a projection $p: \D\to T^*X\oplus J^*A$; note
that $p$ is injective since $\D\cap ((TX\oplus \{0\}) \times
A')=\{0\}$ by  $ii)$ in Prop.\ \ref{prop:equivMP}, and it is onto by dimension count. Now $u_f$
is defined by $p^{-1}(df)$ followed by the projection $K\to TX$,
so it is smooth. Hence $C^\infty(X)^A$ agrees with the algebra of
admissible functions, so it is a Poisson algebra by
Lemma~\ref{lem:pois}.
\end{pf}

\begin{numrmk}
We can fix an isotropic
splitting $j:A^*\to E$, and let $\Pi_X^j$ be the associated
quasi-Poisson bivector field on $X$. The vector field associated
with $f\in C^\infty(X)^A$, constructed in
Prop.~\ref{prop:poisalg}, is $u_f=i_{df}\Pi^j_X$ (it is
\textit{independent} of $j$). The canonical Poisson structure on
$C^\infty(X)^A$ acquires a concrete expression in terms of
$\Pi^j_X$: for $f,g\in C^\infty(X)^A$, $\{f,g\}=\Pi^j_X(df,dg)$.
On the other hand, when $E$ is exact and a splitting $s$ is fixed,
$X$ inherits a Dirac structure $L^s_X$, and $C^\infty(X)^A$ is a
Poisson subalgebra of its algebra of admissible functions (in the
sense of \cite{Co}), characterized by functions $f$ admitting a
Hamiltonian vector field $u_f$ satisfying $dJ(u_f)=0$.
\end{numrmk}

One can also perform moment map reduction for Hamiltonian spaces
for Manin pairs $(E,A)$:

\begin{prop}\label{prop:reduc}
Let $(X,J,\D)$ be a Hamiltonian space for a Manin pair $(E,A)$
over $S$. Let $\mathcal{O}\subset S$ be an orbit of $A$ (i.e., an
integral leaf of the distribution $\rho(A)\subseteq TS$) such that
$J:X\to S$ is transverse to $\mathcal{O}$. Then the $A$-action on
$X$ is tangent to the submanifold $J^{-1}(\mathcal{O})\subset X$,
and the space $C^\infty(J^{-1}(\mathcal{O}))^A$ of $A$-invariant
functions on $J^{-1}(\mathcal{O})$ inherits a Poisson bracket for
which the restriction $C^\infty(X)^A\to
C^\infty(J^{-1}(\mathcal{O}))^A$ is a Poisson map.
\end{prop}

\begin{pf}
From Prop.~\ref{induced-action}, part $i)$, we know that
$dJ(\rho_X(a))=\rho(a)$ for all $a \in J^*A$. This shows that the
$A$-action is tangent to $J^{-1}(\mathcal{O})$. Take $f\in
C^\infty(J^{-1}(\mathcal{O}))^A$, and let $\widetilde{f}$ be an
arbitrary extension of $f$ to $X$. We claim that there exists
$u_{\widetilde{f}} \in \mathfrak{X}(J^{-1}(\mathcal{O}))$
satisfying the condition
\begin{equation}\label{eq:cond}
((u_{\widetilde{f}},d\widetilde{f})_x,0)\in \D,\;\; x\in
J^{-1}(\mathcal{O}).
\end{equation}
(The condition determines $u_{\widetilde{f}}$ uniquely by $i)$  in Prop.\ \ref{prop:equivMP}.)
Indeed, recall from the proof of Prop.~\ref{prop:poisalg} that the
projection $\D\to T^*X\oplus J^*A$ induced by the choice of a
complement $A'$ of $A$ in $E$ is an isomorphism. So one can find
$\widetilde{u}_{\widetilde{f}} \in \mathfrak{X}(X)$ and $a'\in
\Gamma(J^*A')$ such that
$((\widetilde{u}_{\widetilde{f}},d\widetilde{f}),a')\in K$. Since
$f$ is $A$-invariant, we have $a'|_{J^{-1}(\mathcal{O})}=0$ (by
Lemma~\ref{lem:invadm} and $ii)$ in Prop.\ \ref{prop:equivMP}),
which implies that $dJ(\widetilde{u}_{\widetilde{f}})=0$ over
$J^{-1}(\mathcal{O})$. Hence
$u_{\widetilde{f}}:=\widetilde{u}_{\widetilde{f}}|_{J^{-1}(\mathcal{O})}$
is tangent to $J^{-1}(\mathcal{O})$ and satisfies \eqref{eq:cond}.

Given $g \in C^\infty(J^{-1}(\mathcal{O}))^A$, with extension
$\widetilde{g}$ to $X$, then $\D$ being isotropic implies that
$\Lie_{u_{\widetilde{f}}} g = -\Lie_{u_{\widetilde{g}}} f$. As a
result, the bracket $\{f,g\}:=\Lie_{u_{\widetilde{f}}}g =
-\Lie_{u_{\widetilde{g}}} f$ is well-defined, i.e., it only
depends on $f, g$, and not on their extensions. To check the
Jacobi identity, we use $d2)$ in Def.~\ref{def:dirac}. For a given
extension $\widetilde{f}$ of $f\in
C^\infty(J^{-1}(\mathcal{O}))^A$, we can always find a section $e
\in \Gamma(\{0\}\times E)=C^\infty(X\times S, E)$ such that
$\xi_{\widetilde{f}}\in \Gamma(\TX\times E)$, given at $(x,y)\in
X\times S$ by
$((\widetilde{u}_{\widetilde{f}},d\widetilde{f})_x,e_{(x,y)})$,
satisfies $(\xi_{\widetilde{f}})|_{\mathrm{graph}(J)}\in \D$ and
$e_{(x,J(x))}=0$ if $x\in J^{-1}(\mathcal{O})$. From Remark
\ref{rmk:bases}, we see that the restriction of
$\Cour{\xi_{\widetilde{f}},\xi_{\widetilde{g}}}$ to the
submanifold $\{(x,J(x)), \; x\in J^{-1}(\mathcal{O})\}$ gives
$(([u_{\widetilde{f}},u_{\widetilde{g}}],\Lie_{u_{\widetilde{f}}}g),0)
\in \D$, which implies the Jacobi identity. Hence the algebra
$C^\infty(J^{-1}(\mathcal{O}))^A$ has a canonical Poisson bracket,
and it is clear from the construction that the restriction
$C^\infty(X)^A\to C^\infty(J^{-1}(\mathcal{O}))^A$ is a Poisson
map.
\end{pf}

\begin{numrmk}
When an isotropic splitting $j$ of $E$ is
fixed, then the vector field $u_{\widetilde{f}}$ in the proof
Prop.~\ref{prop:reduc} is given by $i_{d\widetilde{f}}\Pi_X^j$,
and the Poisson bracket on $C^\infty(J^{-1}(\mathcal{O}))^A$ can
be computed by
$\{f,g\}:=\Pi_X^j(d\widetilde{f},d\widetilde{g})|_{J^{-1}(\mathcal{O})}$
(one can also directly verify that this is independent of
extensions and that it is a Poisson bracket by \eqref{eq:qpois1}
and \eqref{eq:qpois2}). If $E$ is an exact Courant algebroid, one
can alternatively describe the Poisson algebra
$C^\infty(J^{-1}(\mathcal{O}))^A$ via Dirac reduction, by
identifying this algebra with the admissible functions of the
pull-back image of $L^s_X$ to $J^{-1}(\mathcal{O})$, see
\cite[Sec.~4.4]{BC} (c.f \cite[Sec.~6.4]{BC2}).
\end{numrmk}

\section{Examples from Lie groupoids}\label{sec:integ}

In this section, we consider examples of Manin pairs and
Hamiltonian spaces arising from the theory of Lie groupoids.
Specifically, we let $(E,A)$ be a Manin pair over $S$, and take
$\G\gpd S$ to be the source-simply connected Lie groupoid
integrating $A$, viewed as a Lie algebroid (which we assume to be
integrable). This groupoid may acquire different geometrical
structures via integration: on the one hand, the choice of a
splitting $j$ as in \eqref{eq:j} determines a bivector field $\Pi
\in \mathfrak{X}^2(\G)$ making $\G$ a \textit{quasi-Poisson
groupoid} \cite{ILX} (integrating the Lie quasi-bialgebroid
\eqref{eq:liequasi}); on the other hand, if $E$ is exact and $s$
is a splitting of \eqref{eq:exact}, then $A$ is identified with a
Dirac structure $L_S\subset \TS$ which, according to \cite{BCWZ},
integrates to a 2-form $\omega\in \Omega^2(\G)$ making $\G$ a
\textit{presymplectic groupoid}. The goal of this section is to
establish a direct relationship between the two types of
integration, $\Pi$ and $\omega$, in the light of
Thm.~\ref{thm:equiv}.

\subsection{Multiplicative Dirac structures}\label{subsec:mult}

Let $\hH\gpd \hH_0$ be a Lie groupoid, and let $\G$ be an embedded
subgroupoid of $\hH$, with inclusion homomorphism
$f:\G\hookrightarrow \hH$ (over $f_0:\G_0 \hookrightarrow \hH_0$).
We consider the associated tangent and cotangent Lie groupoids
$T\hH\gpd T\hH_0$ and $T^*\hH\gpd A^*(\hH)$, see e.g.
\cite{MackenzieX:1994}. We observe that there are also natural Lie
groupoid structures on the pull-back bundles $f^*T\hH$, $f^*T^*
\hH$, and on the direct sums $\T\hH=T\hH\oplus T^*\hH$,
$f^*T\hH\oplus f^*T^*\hH$, and $f^*T\hH\oplus T^*\G$.

An (almost) Dirac structure $L$ on $\hH$ supported on $f(\G)$ is
called \textbf{multiplicative} if $L\subset f^*\T\hH$ is a
subgroupoid, i.e., closed under multiplication and inversion. For
$\G=\hH$, this unifies the usual notions of multiplicative
bivector field (when $L=\mathrm{graph}(\Pi)$, for $\Pi\in
\mathfrak{X}^2(\G)$) and multiplicative 2-form (when
$L=\mathrm{graph}(\omega)$, for $\omega\in \Omega^2(\G)$), see
e.g. \cite{MackenzieX:1994}.

Let us consider the map $\mathfrak{F}_f:\mathrm{Dir}(\G)\to
\mathrm{Dir}(\hH)_{f(\G)}$ as defined in
Prop.~\ref{prop:correspondence} (however, in the present
situation, we may consider just almost Dirac structures).

\begin{lem}\label{lem:1}
If $L \subset \TG$ is a multiplicative Dirac structure, then
$\mathfrak{F}_f(L)\subset f^*\T\hH$ is a multiplicative Dirac
structure supported on $f(\G)$.
\end{lem}

\begin{pf} As in the proof of Prop.~\ref{prop:correspondence}, we
consider the maps $\psi:\TG \to f^*T\hH\oplus T^*\G$,
$\psi(u,\alpha)=(df(u),\alpha)$, and $p:f^*\T\hH \to f^*T\hH\oplus
T^*\G$, $p(v,\beta)=(v,f^*\beta)$. Since the maps $df:T\G\to
f^*T\hH$ and $df^*:f^*T^*\hH\to T^*\G$ are groupoid morphisms, so
are $\psi$ and $p$. Hence $\mathfrak{F}_f(L)=
p^{-1}(\psi(L))\subset f^*\T\hH$ is a subgroupoid.
\end{pf}

Let us consider Lie groupoids $\G \gpd \G_0$, $\rR\gpd \rR_0$, and
a homomorphism $J:\G \to \rR$. Then $f:\G\hookrightarrow \G\times
\rR$, $g\mapsto (g,J(g))$ defines an embedded Lie subgroupoid. Let
$\D$ be a multiplicative (almost) Dirac structure on
$\hH:=\G\times \rR$ supported on $f(\G)=\mathrm{graph}(J)$, i.e.,
$\D\subset \TG\oplus J^*\T \rR$ is a subgroupoid. Let $C_\rR
\subset \T \rR$ be a multiplicative almost Dirac structure on
$\rR$.

\begin{lem}\label{lem:2}
Assume that $\D$ and $C_\rR$ satisfy the transversality condition
$(0\oplus J^*C_\rR)\cap \D =\{0\}$ in $\TG\oplus J^* \T\rR$. Then
the composition $\D\circ C_\rR$ is a multiplicative almost Dirac
structure on $\G$.
\end{lem}

\begin{pf}
Let $\pr_\G: \TG\oplus J^*\T\rR \to \TG$ denote the natural
projection. The composition $\D\circ C_\rR$,
$$
(\D\circ C_\rR)_g:=\{(u,\alpha) \in (\TG)_g\;|\; \exists\,
(v,\beta)\in (C_\rR)_{J(g)} \; \mathrm{ s.t. }\;
((u,\alpha),(v,\beta))\in \D_{(g,J(g))} \},\;\; g\in \G,
$$
can be written as $\pr_\G((\TG\oplus J^*C_\rR)\cap \D)$. Since
both $\TG\oplus J^*C_\rR$ and $\D$ are subgroupoids of $\TG\oplus
J^*\T\rR$, so is the intersection $(\TG\oplus J^*C_\rR)\cap \D$.
The transversality condition $(0\oplus J^*C_\rR)\cap \D =\{0\}$ in
$\TG\oplus J^* \T\rR$ says that the restriction of $\pr_\G$ to
$(\TG\oplus J^*C_\rR)\cap \D$ is an \textit{isomorphism} onto
$\D\circ C_\rR$. Since dim($(\D\circ C_\rR)_g$)=dim($\G$) at all
$g$, it follows that $(\TG\oplus J^*C_\rR)\cap \D$ has constant
rank, so it is a smooth vector bundle, and hence $\D\circ C_\rR$
is a smooth vector bundle. Since $\pr_\G$ is a groupoid morphism,
$\D\circ C_\rR$ is a subgroupoid of $\TG$.
\end{pf}

\subsection{Presymplectic and quasi-Poisson groupoids}

We now resume the discussion about presymplectic and quasi-Poisson
groupoids. Our set-up is a Manin pair over $S$ for which the
Courant algebroid is exact. We fix splittings $j$ and $s$ as in
Section \ref{subsec:twoclasses}, so that, after identifications,
we have the following situation: $S$ is equipped with a Dirac
structure $L_S$ (integrable with respect to $\phi_S$), $C_S$ is an
almost Dirac structure such that $\TS=L_S\oplus C_S$. Let us
assume that $L_S$, viewed as a Lie algebroid, admits an
integration to a Lie groupoid, and let $\G\gpd S$ be the
source-simply-connected Lie groupoid integrating $L_S$, with
source/target maps denoted by $\s, \t$, and inversion $i:\G\to
\G$. We know from \cite{BCWZ} that, since $L_S\subset \TS$ is a
Dirac structure, $\G$ has a 2-form $\omega$ making it into a
\textit{presymplectic groupoid}; on the other hand, since $L_S$,
$C_S$ define a Lie quasi-biagebroid, $\G$ inherits a bivector
field $\Pi$ making it a \textit{quasi-Poisson groupoid}
\cite{ILX}.

It results from $(\G,\omega)$ being a presymplectic groupoid that
$J=(\t,\s):\G\to S\times S^{\mathrm{op}}$ is a strong Dirac map
\cite{BCWZ,Xu}, where $S^{\mathrm{op}}$ is equipped with the Dirac
structure $L_S^{\mathrm{op}}:=\{(v,\beta)\;|\; (v,-\beta) \in L_S
\}$ (integrable with respect to $-\phi_S$). In other words,
$(\G,\mathrm{graph}(\omega),J)$ is a Hamiltonian space for the
Dirac manifold $S\times S^{\mathrm{op}}$.
The induced action $\rho_\G$ of $A=L_S\times L_S^{\mathrm{op}}$ on
$J:\G\to S\times S^{\mathrm{op}}$ is
\begin{equation}\label{eq:action}
\rho_\G(u,v)=r_g(u)-l_g(v), \;\;\; u\in (L_S)_{\t(g)},\; v\in
(L_S^\mathrm{op})_{\s(g)},
\end{equation}
where $r_g$ and $l_g$ denote right/left translations on $\G$ (as
well as their tangent maps).

By Thm. \ref{thm:equiv}, there is a bivector field $\Pi_\G\in
\mathfrak{X}^2(\G)$ corresponding to $\omega$, and making $\G$
into a Hamiltonian quasi-Poisson space for the Lie
quasi-bialgebroid determined by the pair $A=L_S\times
L^{\mathrm{op}}_S, A^*\cong C_S\oplus C_S^{\mathrm{op}}$ (with
action \eqref{eq:action} and moment map $J$).

\begin{them}\label{thm:prequasi} The quasi-Poisson bivector field $\Pi_\G$
corresponding to $\omega$ via Thm.~\ref{thm:equiv} agrees with the
multiplicative bivector field $\Pi$ integrating the Lie
quasi-bialgebroid determined by the splitting $\TS=L_S\oplus C_S$.
\end{them}

\begin{pf}
We denote by $\chi \in \Gamma(\wedge^3 L_S)$ the 3-tensor and by
$d_*$ the quasi-differential on $\Gamma(\wedge L_S)$ defined by
the Lie quasi-bialgebroid $\TS=L_S\oplus C_S$. According to
\cite[Thm.~2.34]{ILX}, in order to prove that $\Pi_\G=\Pi$, one
must check that $\Pi_\G$ is multiplicative, and that it satisfies
\begin{equation}\label{eq:comp2}
(d_*f)^r=-[\Pi_\G,\t^*f],\;\;\; \mbox{ and }\;\;
(d_*a)^r=-[\Pi_\G,a^r],
\end{equation}
for all $f\in C^\infty(S)$ and $a\in \Gamma(L_S)$ (for $\xi\in
\wedge^k(L_S)$, $\xi^r$ is defined by
$(\xi^r)_g=r_g(\xi_{\t(g)})$, $g\in \G$).

To show that $\Pi_\G$ is multiplicative, note first that
$J=(\t,\s):\G\to S\times S$ is a groupoid morphism, where $\rR=
S\times S$ is viewed as the pair groupoid. Hence
$\widehat{J}=(\mathrm{Id},J):\G \to \hH$, where $\hH:=\G\times
\rR$, is an embedding which is a groupoid morphism. If follows
from Lemma~\ref{lem:1} that $\D:=\mathfrak{F}_{\widehat{J}}(L_\G)$
is a multiplicative Dirac structure in $\T\hH=\TG\times \T\rR$
supported on $\mathrm{graph}(J)$. Taking $C_\rR= C_S\times
C_S^{\mathrm{op}}$ (writing out the groupoid structure on $\T\rR$,
for $\rR=S\times S$ the pair groupoid, one directly sees that
$C_\rR$ is a multiplicative almost Dirac structure), it follows
from Prop.~\ref{induced-action}, part $iii)$, and
Lemma~\ref{lem:2} that the composition $K\circ C_\rR \subset \TG$
is multiplicative. This implies that $\Pi_\G$ is multiplicative by
\eqref{eq:piX2}. On the other hand, the conditions in
\eqref{eq:comp2} follow directly from \eqref{eq:qpois2} and
\eqref{eq:qpois3} (using the action \eqref{eq:action}).
\end{pf}

One can also verify, using \eqref{eq:qpois1}, that
$\frac{1}{2}[\Pi_\G,\Pi_\G]= \chi^r-\chi^l,$
where $(\chi^l)_g=l_g(i(\chi_{\s(g)}))$, in accordance with
\cite{ILX}.

Thm.~\ref{thm:equiv} provides an explicit construction going from
presymplectic to quasi-Poisson groupoids (and vice-versa), thereby
relating the two integration problems. In the particular case of a
twisted symplectic groupoid (i.e., when $\omega$ is
nondegenerate), the Dirac structure $L_S$ must be the graph of a
bivector field \cite{CaXu} (i.e., it is a twisted Poisson
structure in the sense of \cite{SeWe01}), so we can take $C_S$ to
be $TS$. In this case, the quasi-Poisson bivector field
constructed in Thm.~\ref{thm:equiv} is just the inverse of
$\omega$, and Thm.~\ref{thm:prequasi} recovers
\cite[Prop.~4.5]{ILX}.


\begin{footnotesize}

\end{footnotesize}

\begin{thebibliography}{99}

\bibitem{ABM} A. Alekseev, H. Bursztyn and E. Meinrenken, Pure spinors on Lie groups,
arXiv:0709.1452 [Math.DG].

\bibitem{AK} A. Alekseev and Y. Kosmann-Schwarzbach, Manin Pairs and
moment maps, {\em J. Differential Geom.} {\bf 56} (2000) 133-165.

\bibitem{AKM} A. Alekseev, Y. Kosmann-Schwarzbach and E. Meinrenken,
Quasi-Poisson manifolds, {\em Canad. J. Math.} {\bf 54} (2002) 3-29.

\bibitem {AMM}
{ A. Alekseev, A. Malkin and E. Meinrenken, }\newblock {Lie group
valued moment maps}.
\newblock {\em J. Differential Geom. } {\bf 48} (1998), 445--495.

\bibitem{AX} A. Alekseev and P. Xu, Derived brackets and Courant
algebroids, {\em unpublished manuscript}.

\bibitem{BCMZ} H. Bursztyn, A. Cattaneo, R. Mehta, M. Zambon, in
preparation.

\bibitem{BC}  H. Bursztyn and M. Crainic, Dirac structures,
momentum maps and quasi-Poisson manifolds, {\em The breadth of
symplectic and Poisson geometry} 1-40, Progr. Math. \textbf{232},
Birkh\"auser Boston, 2005.

\bibitem{BC2} H. Bursztyn and M. Crainic, Dirac geometry,
quasi-Poisson actions, and $D/G$-valued moment maps.
Arxiv:0710.0639 [Math.DG].

\bibitem{BCS} H. Bursztyn, M.
Crainic and P. \v Severa, Quasi-Poisson structures as Dirac
structures, {\em Travaux Math\'ematiques} XVI (2005), 41-52.

\bibitem {BCWZ}
H. Bursztyn, M. Crainic, A. Weinstein and C. Zhu, Integration of
twisted Dirac brackets. {\em Duke Math. J.} {\bf 123} (2004),
549--607.

\bibitem{BR}
H. Bursztyn and O. Radko, Gauge equivalence of Dirac structures
and symplectic groupoids. {\em Ann. Inst. Fourier (Grenoble)} {\bf
53} (2003), 309--337.

\bibitem{CaXu}
A. Cattaneo and P. Xu, Integration of twisted Poisson structures.
{\em J. Geom. Phys.} {\bf 49} (2004), 187--196.


\bibitem{Co}
T. Courant, Dirac manifolds, {\em Trans. A.M.S.} {\bf 319} (1990)
631-661.


\bibitem{ILX}  D. Iglesias Ponte, C.
Laurent-Gengoux and P. Xu, Universal lifting theorem and
quasi-Poisson groupoids, {\sl preprint (2005), math.DG/0507396}.

\bibitem{IX}
D. Iglesias Ponte and P. Xu, in preparation.


\bibitem {LWX}
{Z.-J. Liu, A. Weinstein and P. Xu,}\newblock {Manin triples for
Lie algebroids}.
\newblock {\em J. Differential Geom.} {\bf 45}, (1997), 547--574.


\bibitem{Lu}
{J.-H. Lu, Momentum maps and reduction of Poisson actions.} In:
{\em Symplectic geometry, groupoids and integrable systems
(Berkeley, CA, 1989)}, 291-311. Springer, New York, 1991.

\bibitem{Mk}
K.~Mackenzie,
\newblock{\em General theory of Lie groupoids and
Lie algebroids.} London Math. Soc. Lecture Notes Series {\bf 213}.
\newblock{Cambridge University Press, Cambridge}, 2005.

\bibitem{MackenzieX:1994}
K.~Mackenzie and P. Xu, Lie bialgebroids and {P}oisson groupoids,
{\em Duke Math.~J.} {\bf 73} (1994) 415--452.

\bibitem{Roy0}
D. Roytenberg, Courant algebroids, derived brackets and even
symplectic manifolds, {\em U. C. Berkeley thesis}, 1999.
Arxiv:math.DG/9910078.

\bibitem{Roy}
D. Roytenberg, Quasi-Lie bialgebroids and twisted Poisson
manifolds, {\em Lett. Math. Phys.} {\bf 61} (2002) 123--137.

\bibitem{Roy2}
D. Roytenberg, On the structure of graded symplectic
supermanifolds and Courant algebroids, In: Quantization, Poisson
Brackets and Beyond, Theodore Voronov (ed.), {\em  Contemp.
Math.}, {\bf 315}, Amer. Math. Soc., Providence, RI, 2002.

\bibitem{Se}
{ P. \v{S}evera,}\newblock{Letters to A. Weinstein}, available at
http://sophia.dtp.fmph.uniba.sk/$\sim$severa/letters/.

\bibitem{Se2}
{ P. \v{S}evera,} Noncommutative differential forms and quantization of
the odd symplectic category, \emph{Lett. Math. Phys.} 68 (2004), no. 1, pp.
31-39, 2004.

\bibitem{Some}
{ P. \v{S}evera,} Some title containing the words ``homotopy'' and
``symplectic'', e.g.\ this one, \emph{Travaux Math\'ematiques} XVI
(2005), arXiv:math/0105080.

\bibitem {SeWe01}
{P. \v{S}evera and A. Weinstein, }\newblock {Poisson geometry with
a $3$-form background}.
\newblock {\em Prog. Theo. Phys. Suppl.}  {\bf 144} (2001), 145--154.

\bibitem{Xu}
{P. Xu, Morita equivalence and momentum maps}, {\em J.
Differential Geom.} {\bf 67} (2004) 289--333.

\end{thebibliography}
\end{document}